\documentclass[preprint,review,12pt]{elsarticle}

\usepackage[utf8]{inputenc}
\usepackage[T1]{fontenc}
\usepackage{amsmath}
\usepackage{amssymb}
\usepackage{amsfonts}
\usepackage {graphicx}
\usepackage {epstopdf}
\usepackage {subfigure}
\usepackage{epsfig}
\usepackage{color}
\usepackage{mathptmx}
\usepackage{sectsty}
\usepackage{nicefrac}
\usepackage{enumerate}
\usepackage[driverfallback=dvipdfm]{hyperref}
\makeatletter
\def\ps@pprintTitle{%
	\let\@oddhead\@empty
	\let\@evenhead\@empty
	\let\@oddfoot\@empty
	\let\@evenfoot\@oddfoot
}
\makeatother

\begin{document}
	%*****************************
\begin{frontmatter}
\title{Numerical Solution of the Savage–Hutter Equations for Granular Avalanche Flow using the Discontinuous Galerkin Method}
\author[a]{Abdullah Shah \corref{cor1}}
\ead{abdullah\_shah@comsats.edu.pk}
\author[b]{M. Naveed Zafar}
\author[c]{Yulong Du}
\author[d,e]{Li Yuan }
%\ead{lyuan@lsec.cc.ac.cn}
\cortext[cor1]{Corresponding author}
\address[a]{Department of Mathematics, COMSATS University Islamabad, Park Road, Islamabad-45550, Pakistan.}
\address[b]{Earth System Physics, International Centre for Theoretical Physics (ICTP), Trieste 34151, Italy.}
\address[c]{School of Mathematical Sciences, Beihang University, Beijing 100191, P. R. China.}
\address[d]{ICMSEC and LSEC, Academy of Mathematics and Systems Science, Chinese Academy of Sciences, Beijing 100190, P. R. China.}
\address[e]{School of {Mathematical Sciences}, University of Chinese Academy of Sciences, Beijing 100190, P. R. China.}

\begin{abstract}

The Savage-Hutter (SH) equations are a hyperbolic system of nonlinear partial differential equations describing the temporal evolution of the depth and depth-averaged velocity for modelling the avalanche of a shallow layer of granular materials on an inclined surface. These equations admit the occurrence of shock waves and vacuum fronts as in the shallow-water equations while possessing the special reposing state of granular material. In this paper, we develop a third-order Runge-Kutta discontinuous Galerkin (RKDG) method for the numerical solution of the one-dimensional SH equations. We adopt a TVD slope limiter to suppress numerical oscillations near discontinuities. And we give numerical treatments for the avalanche front and for the bed friction to achieve the well-balanced reposing property of granular materials.  Numerical results of the avalanche of cohesionless dry granular materials down an inclined and smoothly transitioned to horizontal plane under various internal and bed friction angles and slope angles are given to show the performance of the present numerical scheme.
\end{abstract}

\end{frontmatter}

\section{Introduction}

Landslides, snow avalanches,  debris flows, and pyroclastic flows are destructive natural phenomena that may cause a massive amount of life and property losses in mountainous regions. For these geophysical mass flows, given the initial piles, the prediction of the flowing velocity, run-out zone, and deposit distribution is of great importance in hazard assessments \cite{Pud2007book}. These phenomena can be mathematically modeled by using a continuum mechanical or discrete mechanical approach. In the continuum mechanical approach for debris flows, the moving mixture of sediment and water can be treated as a continuum fluid, allowing the use of the Navier-Stokes equations. Additionally, these flows commonly exhibit the characteristics of shallowness, i.e., the flow depth normal to the basal topography is relatively small compared with the lateral spreading scale of the avalanche, and the lateral velocities are more significant than the normal velocity. Thus depth-averaged shallow granular flow models \cite{SH89, Gray1999, Pudasaini2003},
which are extensions of the traditional shallow water equations \cite{Vre1994book, Kuba2009}, have been developed to model granular avalanche flows. These models are derived by integrating the incompressible Navier-Stokes equations from the basal surface along the normal direction to the free surface, and they consist of two partial differential equations for the temporal evolution of the depth and depth-averaged velocities tangential to the bed (cf. \cite{Fei2015}).  Initially, Savage and Hutter proposed to use the Mohr-Coulomb soil constitutive law for dry granular mass in their seminal Savage-Hutter (SH) model \cite{SH89}.
However, it is recognized that real avalanches are multiphase flows with fluids and various sizes of solid grains.  To account for fluid effects, depth-averaged mixture fluid models \cite{Iverson2001} and two-fluid models
\cite{Pitman} were introduced. Meanwhile,  developments of the SH models have been carried out in different directions, including formulations in two-dimensional curved and twisted channel topography \cite{Hutter1993, Gray1999,Pudasaini2003} and arbitrary topography \cite{Bouchut2004,Luca1,Luca2,LiYuan2018},  multi-layer flow models \cite{Taka2005}, basal erosion/deposit models \cite{tai2008}, and GIS-based parallel adaptive computation \cite{Patra05}, to name a few.

Since the hyperbolic properties of the SH equations are similar to those of the shallow water equations, many numerical methods established for the latter could applied to solve the SH equations. Significant efforts have been made in the past three decades. Earlier simulation studies used Lagrangian methods \cite{Gray1999,Meilong}.
Later, Wang et al. \cite{wang2004} used a high-resolution Non-Oscillatory Central (NOC) difference scheme with a 2nd-order MUSCL or 3rd-order WENO reconstruction. A finite volume scheme with a Roe type approximate Riemann solver was developed by  Pelanti et al. \cite{Pelanti-ESAIM08} for  one-dimensional two-phase shallow granular flows. A Godunov type finite volume scheme was used by  Xia and Liang \cite{Xilin18}. A finite volume scheme with the MUSCL reconstruction and Harten-Lax-van Leer Contact (HLLC) numerical flux was used by Zhai et al. \cite{liyuan15} for dry shallow granular flows. Listed above are only  a few of numerous existing finite volume methods for the SH models. Also, gas kinetic schemes were used by Mangeney et al. \cite{Mangeney2003} for modelling dry granular avalanches, and by Chen et al. \cite{GKS2013} for the anisotropic SH equations.

While the SH equations are similar to the shallow water equations in the sense that both systems allow shock waves and wet-dry fronts to occur, there are some differences due to the solid-like constitutive law used in the SH equations. First, the lateral motion is not isotropic due to different earth pressure coefficients.  Second, the granular materials can keep/regain a reposing state with an inclined free surface as long as the inclination angle of the free surface is less than the internal friction angle. Numerical schemes must take this solid-like behavior into account in the discretizations of the momentum equations \cite{Mangeney2003,liyuan15}.

In recent years, high-order discontinuous Galerkin (DG) finite element methods \cite{DGbook2000} have gained much attention in several fields such as computational fluid dynamics \cite{PeterBastian18}, computational acoustics \cite{Schoeder19}, computational electromagnetics \cite{Hongxin19}, etc.  This is spurred by the introduction of the RKDG methods for hyperbolic conservation laws by Cockburn and Shu \cite{CS98}. RKDG methods use piecewise polynomials to approximate the solution, and after the spatial DG discretization in weak form, a system of ordinary differential equations is obtained and solved using a high-order strong stability-preserving (SSP) Runge-Kutta method \cite{Shu88,Shu-Osher}. RKDG methods have many advantages such as local conservation, higher-order accuracy, compactness,  easy imposition of boundary conditions,  suitability for complex meshes, parallelization and adaptivity, and provable convergence (cf. review article \cite{cockshureview01}).
RKDG methods have been widely used for the shallow water equations \cite{Kuba2009, Kess10, Posi10,Khanbook14, Li2017}. %[{\color{red}\cite{qamar2014}is not DG, replace}]

In this work, we apply a 2nd-degree polynomial RKDG method for the numerical solution of the one-dimensional SH equations \cite{Pud2007book} to study the spreading of a granular avalanche down an inclined plane \cite{Hutter1993}. The time integration is done by the 3rd-order SSP Runge-Kutta method. Special attention is paid to the discrete treatment of the static equilibrium state of Coulomb materials.

The organization of this paper is as follows.  In section 2, the governing equations of a 1D SH model are given. In Section 3,  numerical discretizations both in space and time are given together with
the slope limiter and numerical treatments of dry-wet fronts and reposing states of Mohr-Coulomb type plastic materials. The  computed results of some test cases are provided in Section 4 to validate the proposed scheme and illustrate the effects of the phenomenological parameters and inclination angle on the flow behavior. Section 5 concludes this work.

\section{Savage-Hutter shallow granular flow model}
The 1D dimensionless  SH equations for the shallow flow of a finite mass of granular material down an inclined plane can be written as \cite{tai2001}
\begin{subequations}
	\label{1b}
	\begin{eqnarray}
		\frac{\partial h}{\partial t}+\frac{\partial q}{\partial x} &=&0,
		\label{eq1} \\
		\frac{\partial q}{\partial t}+\frac{\partial }{\partial x}\left( \frac{q^{2}%
		}{h}+\frac{\beta h^{2}}{2}\right)  &=&hs(u),  \label{eq2}
	\end{eqnarray}
\end{subequations}
where $h$ is the flow thickness and $q=hu$ is the discharge with $u$ is the depth-averaged velocity component in the $x$-direction. The factor $\beta$ is defined as
\begin{equation*}
\beta =\varepsilon \cos \left( \zeta \left( x\right)  \right) K,
\end{equation*}
where $\varepsilon=H/L \ll  1$ is the aspect ratio of the characteristic thickness $H$ to longitudinal extent $L$, $\zeta \left(x\right)$ is the inclination angle of the reference surface \cite{wang2004}, and $K$ is the   earth pressure coefficient defined by the Mohr-Coulomb criterion as follows,
\begin{equation*}
	K:= \frac{p_{xx}}{p_{zz}}=\left\{
	\begin{array}{l}
		K_{\mathrm{act}}=2\sec ^{2}\varphi \left( 1-\sqrt{1-\cos ^{2}\varphi \sec^{2}\delta }\right) -1, \\
		K_{\mathrm{pass}}=2\sec ^{2}\varphi \left( 1+\sqrt{1-\cos ^{2}\varphi \sec^{2}\delta }\right) -1.%
	\end{array}%
	\right.
\end{equation*}

Here, $\varphi$ is the internal friction angle of the granular material and $\delta $ is the friction angle between the avalanche and the base. The subscripts "act" ($-$ sign) and "pass" ($+$ sign) denote  active and passive stress states respectively, which become effective when the avalanche extends or contracts in  the down-slope direction:
\begin{equation*}
K=\left\{ \begin{split}   K_{\rm act}, &~~ \displaystyle  \frac{\partial u}{\partial x}\geq 0, \\
K_{\rm pass}, & ~~ \displaystyle  \frac{\partial u}{\partial x}<0.
\end{split}
\right.
\end{equation*}
The source term $s(u)$ in the momentum equation \eqref{eq2} represents the net driving
acceleration in the down-slope direction, i.e.,
\begin{equation}
	s(u)=\sin \zeta \left( x\right) -\frac{u}{\sqrt{u^{2}}}\tan \delta \cdot {%
		\max }\left( 0,~\cos \zeta \left( x\right) +\chi \kappa u^{2}\right),
	\label{srce}
\end{equation}
where $\kappa =-\dfrac{\partial \zeta \left( x\right) }{\partial x}$ is the
local curvature of the reference surface, and $\chi \kappa $ is local
stretching of the curvature. The first term in Eq.~(\ref{srce}) is the tangential component of the gravity, and the second term is the basal friction.

%If a moving fluid $|u|<10^{-14}$, then $s(u)$ is replaced by the driving force due to the gravity $+$ inertia $+$ surface slope $f_{R}$, and if a motionless fluid $u=0$, then $s(u)$ is replaced by the net driving force due to  the gravity $+$ surface slope $f_{D}$.

In this work, we take $\chi=1$ and the inclination angle $\zeta(x)$ is defined as in \cite{wang2004}; 
\begin{equation}\label{incl}
\zeta \left( x\right) =\left\{
\begin{array}{ll}
\zeta _{0} &    0 \leq x  \leq   17.5,    \\
\zeta _{0}\left(\displaystyle 1-\frac{x-17.5}{4} \right)  &    17.5 < x < 21.5,  \\
0  &  x  \ge 21.5,
\end{array}
\right.
\end{equation}
where $\zeta_0$ will be assigned different values in the  numerical tests.

\begin{figure}[htp]
	\centering
	\includegraphics[width=0.9\textwidth]{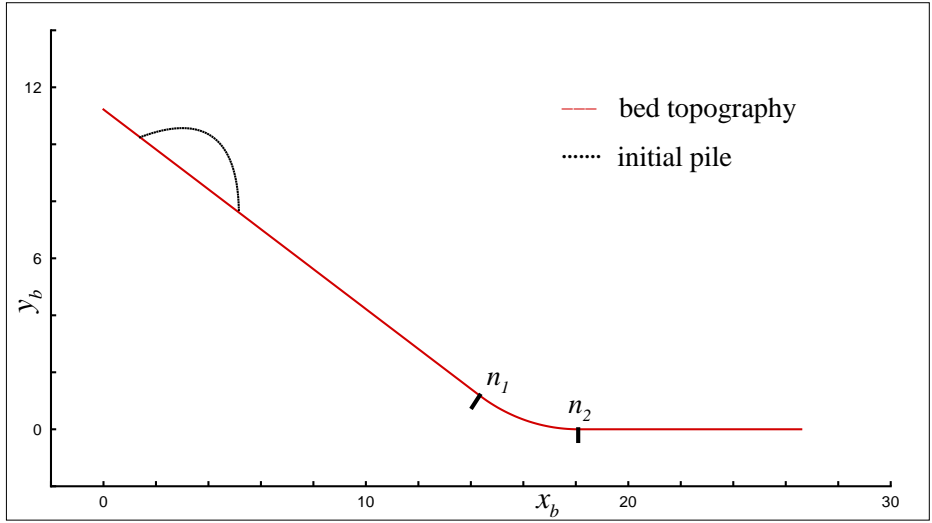}
	\caption {Sketch of the bed topography defined Eq.~(\ref{bed}) with $y_{b}^{0}=10.0$ and initial granular pile. The bed consists of three parts: the inclined part $[0,n_1)$, the transition part $[n_1, n_2]$, and the horizontal run-out part $(n_2,30]$. The transformation $x \to (x_b,y_b)$ is given in the appendix. }\label{top}
\end{figure}

Figure~\ref{top} shows a sketch of  the basal topography defined as the elevation above the reference surface as follows,
\begin{equation} \label{bed}
y_{b}\left(x\right)=y_{b}^{0}\left(1-\cos\left(\frac{\pi}{4}\right) \right) \sin \zeta \left(x\right).
\end{equation}

%\begin{figure}[htp]
%	\centering
%	\includegraphics[width=0.9\textwidth]{xi.eps}
%	\caption {Sketch of the reference surface defined by Eq.~(\ref{incl}) with $\zeta_0=35^{\circ}$.}\label{xi}
%\end{figure}

%\begin{figure}[htp]
%	\centering
%	\includegraphics[width=0.9\textwidth]{xi.eps}
%	\caption {Sketch of the reference surface defined by Eq.~(\ref{incl}) with $\zeta_0=35^{\circ}$.}\label{xi}
%\end{figure}

By introducing the vector of conservative variables $\textbf{U}=({U_1, U_2})^T=(h,q)^T$, the convective flux vector $\textbf{F}(\textbf{U})=({F_1(\textbf{U}),F_2(\textbf{U})})^T=(q,q^{2}/h+\beta
h^{2}/2)^T$, and the source vector $\textbf{S}(\textbf{U})=\left(S_1(\textbf{U}),S_2(\textbf{U})\right)^T=\left(0,hs(u)\right)^T$, the mass and momentum balance equations (\ref{eq1}) and (\ref{eq2}) can be written in the vector form as follows,
\begin{equation}
\frac{\partial \textbf{U}}{\partial t}+\frac{\partial \textbf{F}\left( \textbf{U}\right) }{\partial x}=\textbf{S}\left( \textbf{U}\right), \label{eq3}
\end{equation}
where the Jacobian matrix $A(\textbf{U})=\textbf{F}^{\prime}(\textbf{U})$ is given as
\begin{equation*}
A(\textbf{U})=\frac{\partial \textbf{F}(\textbf{U})}{\partial \textbf{U}}=\left(\begin{array}{cc}
0 & 1 \\ \beta h-u^{2} & 2u
\end{array}\right).
\end{equation*}
The two distinct eigenvalues of $A(\mathbf{U})$ are
\begin{equation*}
	\lambda _{1}(\mathbf{U})=u-\sqrt{\beta h}~,\text{ \ and }\quad \lambda _{2}(\mathbf{U})=u+\sqrt{\beta h}~.
\end{equation*}

%\textbf{(Computational Difficulty: when $h$ approach zero, see Liang QH, Kurganov SWE papers how to compute $u$)}

Note that when the velocity of the fluid  is smaller than the speed $\sqrt{\beta h}$ of the gravity wave, i.e., $|u| < \sqrt{\beta h}$, the flow is said to to be sub-critical and then $\lambda_{1}<0$ and $\lambda_{2}>0$. Under the sub-critical condition, there are no shocks. When $|u| > \sqrt{\beta h}$, the flow is said to be super-critical. Any transition from a supercritical to a sub-critical flow state may produce a shock wave.

 %The ratio $Fr=u/ \sqrt{\beta h}$ is called \textit{Froude number}}.

The right and left eigenvector matrices of $A\mathbf{(U)}$ are respectively
\begin{equation*}
	R=\left(
	\begin{array}{cc}
		1 & 1 \\
		u-\sqrt{\beta h} & u+\sqrt{\beta h}%
	\end{array}%
	\right) ,~R^{-1}=\frac{1}{2\sqrt{\beta h}}\left(
	\begin{array}{cc}
		u+\sqrt{\beta h} & -1 \\
		-u+\sqrt{\beta h} & 1%
	\end{array}%
	\right) .
\end{equation*}
The Jacobian matrix can be diagonalized as $R^{-1}AR=\Lambda $, where $\Lambda=\text{diag}(\lambda_1, \lambda_2) $. As long as $h\geq 0$, the system is hyperbolic.
\section{Runge-Kutta discontinuous Galerkin method}
In this section, we present the RKDG method  for  Eq.~(\ref{eq3}) including numerical treatments of the avalanche front and the basal friction to achieve the reposing state of granular materials.

\subsection{Spatial discretization}

To begin, the one-dimensional computational domain $\Omega=[0, L]$ is divided into $N$ cells with the cell interfaces  $0=x_{1/2}<x_{3/2}< \dots < x_{N+1/2}=L$. The point $x_{j}=\frac{1}{2}(x_{j-1/2}+x_{j+1/2})$ is the center of the $j$-th cell $I_{j}=[x_{j-1/2},x_{j+1/2}]$, and we denote cell size by $\Delta_{j}=x_{j+1/2}-x_{j-1/2}$ with  $\Delta x=\underset{1\leq j\leq N}{\max }\Delta_{j}$.
The function space of piecewise polynomials is given by
\begin{equation*}
	V_{h}^{k}=\left\{ v_h\in L^\infty(\Omega): v_h |_{I_{j}}\in P^{k}\left( I_{j}\right)
	,~~1\leq j\leq N\right\} ,
\end{equation*}
where $ P^{k}(I_{j})$ denotes the polynomial  space of degree at most $k$ on cell $I_{j}$. The functions in $V_{h}^{k}$ are allowed to have discontinuity across cell interfaces. In this study, we choose the scaled
Legendre polynomials $\phi_l(x)$,~$l=0,1,\cdots, k$ as the local basis functions $v_h$ over $I_{j}$, which are local orthogonal basis over $I_j$. This gives the  local basis function space $P^{k}(I_j)$:
\begin{equation*} \label{localspace}
P^{k}(I_j)=\left\{ 1,2\left( \frac{x-x_{j}}{\Delta_{j}}\right) ,\frac{1}{2}
\left( \frac{12\left( x-x_{j}\right) ^{2}}{\Delta_{j}^{2}}-1\right), \cdots, \phi_k(x) \right \}.
\end{equation*}
The  approximate solution $\textbf{U}_h(x,t)=( h_j(x,t),q_j(x,t))^{T} \in  V^k_h$ to the exact $\textbf{U}(x,t)$ over each cell $I_{j}$ can be written as
\begin{equation}
\textbf{U}_h(x,t)=\sum_{l=0}^{k}\textbf{U}_{l}^{j}(t)\phi_l(x)~~~~ \text{for $ x \in I_{j} $}. \label{app}
\end{equation}

Next, we multiply Eq.~(\ref{eq3}) with the test function $\phi_m$, use Eq.~(\ref{app}) and integrate  by parts over the interval $I_{j}$, which results in the following weak form:
\begin{eqnarray}
&&\int_{I_{j}}\partial_t\left(\sum_{l=0}^{k}\textbf{U}^{j}_{l}(t)\phi_{l}(x)\right)
\phi_{m}(x)\mbox{d}x-\int_{I_{j}}\textbf{F}\left( \textbf{U}_h(x,t)\right) \partial_x \phi_{m}(x)\mbox{d}x \nonumber \\
&&+ \hat{\textbf{F}}_{j+\frac 12} \phi_{m}\left(x^-_{j+\frac 12}\right)
 - \hat{\textbf{F}}_{j-\frac 12}  \phi_{m}\left(x^+_{j-\frac 12}\right)\nonumber \\
&& = \int_{I_{j}}\textbf{S}\left( \textbf{U}_h(x,t)\right) \phi_{m}(x) \mbox{d}x, ~~m=0,1,\ldots ,k, \label{weakform}
\end{eqnarray}
with the initial values of the degrees of freedom (DOFs) obtained by the projection
\begin{equation}\label{initia}
\int_{I_{j}}\left(\sum_{l=0}^{k}\textbf{U}^{j}_{l}(0)\phi_{l}(x)\right) \phi_m(x) \mbox{d}x = \int_{I_{j}} \textbf{U}_0(x)\phi_m(x) \mbox{d}x, ~~m=0,1,\dots ,k.
\end{equation}
In Eq.~(\ref{weakform}),~  $\hat{\textbf{F}}_{j+1/2}  =\hat{\textbf{F}}\left(\textbf{U}_{j+1/2}^{-},\textbf{U}_{j+1/2}^{+}\right)$ is a monotone numerical flux obtained  from the exact or approximate Riemann solver,  which connects the discontinuous approximate solutions at the cell interface $x_{j+1/2}$. In this work, the local Lax-Friedrichs flux is used as numerical flux i.e.,
\begin{equation}\label{llf}
\hat{\textbf{F}}\left(\textbf{U}_{j+\frac 12}^{-},\textbf{U}_{j+\frac 12}^{+}\right)=\frac{1}{2}\left[
\textbf{F}\left(\textbf{U}_{j+\frac 12}^{-}\right)+\textbf{F}\left(\textbf{U}_{j+\frac 12}^{+}\right)-\alpha
\left(\textbf{U}_{j+\frac 12}^{+}-\textbf{U}_{j+\frac 12}^{-}\right)\right] .
\end{equation}
 where $\alpha$ is an estimate of the largest absolute value of eigenvalue of the Jacobian matrix in the neighborhood of the interface $x_{j+1/2}$, i.e.,
$$\alpha =\max\left(\lambda_{m}(\textbf{U}_{j+1/2}^{-}),\lambda_{m}(\textbf{U}_{j+1/2}^{+})\right),~~ \lambda_{m}=\max(|\lambda_{1}|,|\lambda_{2}|).$$

The last term in Eq.~\eqref{llf} represents numerical dissipation. The first component of this term do not allow the mass conservation equation to reach steady state when  $h$ gradients are not zero even though the flow velocity is zero.  To resolve it,  we multiply the first component of this term with a factor $\eta_{j+1/2}=\max(M^{\rm stop}_{j}, M^{\rm stop}_{j+1})$, where the tag $M^{\rm stop}_{j} =0 $ indicates a resting cell and  $M^{\rm stop}_{j} =1 $  indicates a flowing cell,  see the end of Sec.~\ref{sec:3.4}. This treatment is inspired by a technique in \cite{Mangeney2003}.

Due to the orthogonality of the Legendre polynomials, we obtain a diagonal mass matrix $M$ for the first term in Eq.~(\ref{weakform}), so this equation can be written as:
\begin{eqnarray}
	&&\left. M  \left(
	\begin{array}{c}
		\displaystyle\frac{\mbox{d}\mathbf{U}_{0}^{j}}{\mbox{d}t} \\
		\vdots  \\
		\displaystyle\frac{\mbox{d}\mathbf{U}_{k}^{j}}{\mbox{d}t}%
	\end{array}%
	\right) -\int_{I_{j}}\mathbf{F}\left( \mathbf{U}_{h}(x,t)\right) \left(
	\begin{array}{c}
		\phi _{0}^{\prime }(x) \\
		\vdots  \\
		\phi _{k}^{\prime }(x)%
	\end{array}%
	\right) \mbox{d}x+\hat{\mathbf{F}}_{j+\frac{1}{2}}\left(
	\begin{array}{c}
		1 \\
		\vdots  \\
		1%
	\end{array}%
	\right) \right.   \notag \\
	&&\left. -\hat{\mathbf{F}}_{j-\frac{1}{2}}\left(
	\begin{array}{c}
		(-1)^{0} \\
		\vdots  \\
		(-1)^{k}%
	\end{array}%
	\right) =\int_{I_{j}}\mathbf{S}\left( \mathbf{U}_{h}(x,t)\right) \left(
	\begin{array}{c}
		\phi _{0}(x) \\
		\vdots  \\
		\phi _{k}(x)%
	\end{array}%
	\right) \mbox{d}x,\right.   \label{weak2}
\end{eqnarray}
where we have used the following properties of the Legendre polynomials:
\begin{equation*}
\phi_{m}\left(x_{j+
\frac{1}{2}}^{-}\right)=1, \quad \text{and} \quad \phi_{m}\left(x_{j-\frac{1}{2}}^{+}\right)=(-1)^{m}.
\end{equation*}
For the degree $k=2$ DG approximation, the mass matrix $M$ is given as
\[
\begin{split}
& \left. M=\left[ \int_{I_{j}}\phi _{l}(x)\phi _{m}(x)\mbox{d} x\right] =\Delta _{j}%
\left[ \frac{\delta _{lm}}{2l+1}\right] =\Delta _{j}\left[
\begin{array}{ccc}
1 & 0 & 0 \\
0 & \frac{1}{3} & 0 \\
0 & 0 & \frac{1}{5}%
\end{array}%
\right] \right., \quad  l,m=0,1,2.
\end{split}%
\]
where
\begin{equation*}
\delta _{lm}=\left \{
\begin{array}{ccc}
1 & \text{if} & l=m, \\
0 & \text{if} & l\neq m.%
\end{array}
\right.
\end{equation*}
The derivatives of the Legendre basis functions are given by
\begin{equation*}
	\left\{ \phi _{0}^{\prime }=0,~~\phi _{1}^{\prime }=\frac{2}{\Delta _{j}}%
	,~~\phi _{2}^{\prime }(x)=\frac{12(x-x_{j})}{\Delta _{j}^{2}},~~\cdots
	\right\} .
\end{equation*}
For evaluating the mass flux  $F_1(\textbf{U}_h(x,t))$ component in the second term in Eq.~(\ref{weak2}),  the following integrals are useful:
\begin{equation*}
\begin{split}
&\int_{I_{j}}\phi_l(x) \phi^{\prime }_0 \mbox{d} x=0, \quad l=0,1,2. \\
&\int_{I_{j}}\phi_0(x) \phi^{\prime }_1 \mbox{d} x=2, \quad
\int_{I_{j}}\phi_l(x) \phi^{\prime }_1 \mbox{d} x =0,\quad l >0  . \\
&\int_{I_{j}} \phi_1 (x) \phi^{\prime }_2 (x) \mbox{d} x =2,\quad
\int_{I_{j}}\phi_l(x) \phi^{\prime }_2 (x) \mbox{d} x =0,\quad l \neq 1.
\end{split}
\end{equation*}
Thus, the $F_1(\textbf{U}_h(x,t))$ component can be integrated exactly (assume  $k=2$)
\begin{eqnarray*}
 \int_{I_{j}}{F_{1}}\left( \mathbf{U}_{h}(x,t)\right) \left(
	\begin{array}{l}
		\phi _{0}^{\prime } \\
		\phi _{1}^{\prime } \\
		\phi _{2}^{\prime }(x)%
	\end{array}%
	\right) \mbox{d}x  &=&  \int_{I_{j}}\left( \sum_{l=0}^{2}q_{l}^{j}(t)\phi _{l}(x)\right)
	\left(
	\begin{array}{l}
		\phi _{0}^{\prime } \\
		\phi _{1}^{\prime } \\
		\phi _{2}^{\prime }(x)%
	\end{array}%
	\right) \mbox{d}x     \notag \\
	& =  & \left[
	\begin{array}{ccc}
		0 & 0 & 0 \\
		2 & 0 & 0 \\
		0 & 2 & 0%
	\end{array}%
	\right] \left[
	\begin{array}{c}
		q_{0}^{j}(t) \\
		q_{1}^{j}(t) \\
		q_{2}^{j}(t)%
	\end{array}%
	\right] .
\end{eqnarray*}
{However, the momentum flux  $F_2(\textbf{U}_h(x,t))$ component in the second term in Eq.~(\ref{weak2}) is highly nonlinear and is not easy to integrate analytically. For the $P^2$-DG solution, we use the following three point Gauss quadrature formula \cite{BurdenFairs}:
\begin{equation}
\int_{I_{j}}F_2\left(\textbf{U}_h(x,t)\right)\phi^\prime_{m} (x) \mbox{d}x \approx \Delta_j \sum\limits_{g=1}^{3}\omega
_{g}F_2\left(\textbf{U}_h(x_{g},t)\right) \phi^\prime_{m} (x_g),~~m=0,1,2,  \label{gauss}
\end{equation}
where $(\omega_{1},\omega_{2},\omega_{3}) =\left(\frac {5}{18}, \frac {8}{18}, \frac {5}{18}\right)$ are the Gauss quadrature weights, and $ (x_j-\sqrt{\frac35}\frac{\Delta_j}{2}$, $x_j$, $ x_j+\sqrt{\frac35} \frac{\Delta_j}{2} )$  are the Gauss points \cite{BurdenFairs}.

%\textcolor[rgb]{1.00,0.00,0.00}{\textbf{[1. Well-balance, 2. Stopping criterion, 3. Accuracy. Read Liang and Qmar phd thesis on RKDG]} }

The integration of the source term in Eq.~(\ref{weak2}) needs careful treatment. The first component $S_1(\textbf{U})$ is zero. The second component $S_2(\textbf{U})$ can be computed  by the same  Gauss quadrature formula as \eqref{gauss},
\begin{equation*}
\int_{I_{j}}hs(u)\phi_{m} (x) \mbox{d}x \approx \Delta_j \sum\limits_{g=1}^{3}\omega
_{g} h(x_g,t) s(x_g,t)  \phi_{m} (x_g),~~m=0,1,2.  \label{gauss2}
\end{equation*}
%\textcolor{red}{or by second-order approximation as
%	\begin{equation*}
%	\begin{split}
%	\int_{I_{j}} {h}s(u)\phi_{m}\mbox{d} x\approx & s(\bar{u}_{j})\int_{I_{j}}\left(
%	\sum_{l=0}^{2}h^{j}_{l}(t)\phi_l(x)\right)  \phi_m(x)  \mbox{d} x \\
%	=&s(\bar{u}_{j}) \frac{\Delta_j }{2m+1}  \sum_{l=0}^2 h^j_{l} (t)  \delta_{lm}, ~m=0,1,2..
%	\end{split}
%	\end{equation*}
%}
%(\textcolor{red}{we can neglict this formula}).
%%After discretizing in space by the DG method, Equation (\ref{weak2}) can be finally rewritten in a semi-discrete form as follows;
%%\begin{equation}
%%\Delta_j \frac{\mbox{d} \textbf{U}^j_0}{\mbox{d} t}=-\left(\hat{\textbf{F}}_{j+\frac 12}-\hat{\textbf{F}}_{j-\frac 12}\right)+\Delta_j \sum_{g=1}^{M=3}\omega_g \textbf{S}^j(x_g,t)
%%\end{equation}
%%[{\color{red} write out each ODE, l=0,1,2  for each DOF for complete, detail H}]
%%  \begin{eqnarray}
%%&&\left. \frac{\Delta_{j} } {2l+1} \frac{\mbox{d}\textbf{U}^j_l}{\mbox{d}t}=\left[\textbf{H}_l (\textbf{U}%
%%)-\left( \hat{\textbf{F}}(U_{j+\frac{1}{2}}^{-},U_{j+\frac{1}{2}}^{+})-(-1)^{l} \hat{\textbf{F}}
%%(U_{j-\frac{1}{2}}^{-},U_{j-\frac{1}{2}}^{+})\right) +\Delta_j \sum_{g=1}^3 \textbf{S}_0(\tilde{U})\right]
%%,\right.   \notag \\
%%&&\quad \left. ~~~~~~~~~~~~~~~~~~~~~~~~~~~~~~~~~~~~~~~~~~~~~~~~~~~~~~~~~~~~~~~~~~~~~~~~~~~l=0,1,2 \dots k.\right.   \notag \\
%%&&\left. {}\right.   \label{semidiscrete}
%%\end{eqnarray}

After discretizing in space by the DG method, Eq. (\ref{weak2}) can finally be rewritten in a semi-discrete form as follows,
\begin{equation}\label{semidiscrete}
	\begin{split}
		\Delta _{j}\frac{\mbox{d}\mathbf{U}_{0}^{j}}{\mbox{d}t} &= \mathbf{H}_{0}(%
		\mathbf{U}_{h})-\left( \hat{\mathbf{F}}_{j+\frac{1}{2}}-\hat{\mathbf{F}}_{j-%
			\frac{1}{2}}\right) +\Delta _{j}\sum_{g=1}^{3}\omega _{g}{\mathbf{S}\left(\textbf{U}_h(x_{g},t)\right)}\phi _{0}(x_{g}), \\
		\frac{\Delta _{j}}{3}\frac{\mbox{d}\mathbf{U}_{1}^{j}}{\mbox{d}t} &=
		\mathbf{H}_{1}(\mathbf{U}_{h})-\left( \hat{\mathbf{F}}_{j+\frac{1}{2}}+\hat{%
			\mathbf{F}}_{j-\frac{1}{2}}\right) +\Delta _{j}\sum_{g=1}^{3}\omega _{g}%
		{\mathbf{S}\left(\textbf{U}_h(x_{g},t)\right)}\phi _{1}(x_{g}), \\
		\frac{\Delta _{j}}{5}\frac{\mbox{d}\mathbf{U}_{2}^{j}}{\mbox{d}t} &=
		\mathbf{H}_{2}(\mathbf{U}_{h})-\left( \hat{\mathbf{F}}_{j+\frac{1}{2}}-\hat{%
			\mathbf{F}}_{j-\frac{1}{2}}\right) +\Delta _{j}\sum_{g=1}^{3}\omega _{g}%
		{\mathbf{S}\left(\textbf{U}_h(x_{g},t)\right)}\phi _{2}(x_{g}),
	\end{split}
\end{equation}
with
\begin{equation*}
	\begin{aligned}
		\mathbf{H}_{0}(\mathbf{U}_{h}) &= \int_{I_{j}}\mathbf{F}\left( \mathbf{U}%
		_{h}(x,t)\right) \phi _{0}^{\prime }(x)\mbox{d}x=0, \\
		\mathbf{H}_{1}(\mathbf{U}_{h}) &= \int_{I_{j}}\mathbf{F}\left( \mathbf{U}%
		_{h}(x,t)\right) \phi _{1}^{\prime }(x)\mbox{d}x=\left(
		\begin{array}{c}
			2q_{0}^{j}(t) \\
			 \Delta _{j}\sum\limits_{g=1}^{3}\omega _{g}{F_{2}}\left( \mathbf{U}%
			_{h}(x_{g},t)\right) \phi _{1}^{\prime }(x_{g})%
		\end{array}%
		\right) , \\
		\mathbf{H}_{2}(\mathbf{U}_{h}) &=\int_{I_{j}}\mathbf{F}\left( \mathbf{U}%
		_{h}(x,t)\right) \phi _{2}^{\prime }(x)\mbox{d}x=\left(
		\begin{array}{c}
			2q_{1}^{j}(t) \\
			\Delta _{j}\sum\limits_{g=1}^{3}\omega _{g}{F_{2}}\left( \mathbf{U}%
			_{h}(x_{g},t)\right) \phi _{2}^{\prime }(x_{g})%
		\end{array}%
		\right) .
	\end{aligned}
\end{equation*}

\subsection{Time integration}

We apply the 3rd-order {SSP Runge-Kutta scheme \cite{Shu88, Shu-Osher, cockshureview01} to solve the ODE system of equations~\eqref{semidiscrete} using the initial values \eqref{initia}}, i.e.,
\begin{equation*}
\dfrac{\mbox{d} \textbf{U}_h  } {\text{d} t}=\pounds (\textbf{U}_h),\quad \textbf{U}_h(0)=\mathbb{P}_{V^k_h} \left(\textbf{U}_0(x)\right).
\end{equation*}%

%\textcolor{red}{ The correct option based on system (13) and Eq (7) may be 
%\begin{equation*}
%\dfrac{\ {d\mathbf{U}_{l}^{j}(t)}}{dt}=\pounds (\textbf{U}_l^j(t))
%\end{equation*}%
The 3rd-order  optimal SSP RK  scheme takes the form:
\begin{eqnarray*}
	\mathbf{U}^{\left( 1\right) } &=&\mathbf{U}^{n}+\Delta t\pounds \left(
	\mathbf{U}^{n},t^{n}\right) , \\
	\mathbf{U}^{\left( 2\right) } &=&\frac{3}{4}\mathbf{U}^{n}+\frac{1}{4}\left(
	\mathbf{U}^{\left( 1\right) }+\Delta t\pounds \left( \mathbf{U}^{\left(
		1\right) },t^{n}+\Delta t\right) \right) , \\
	\mathbf{U}^{n+1} &=&\frac{1}{3}\mathbf{U}^{n}+\frac{2}{3}%
	\left( \mathbf{U}^{\left( 2\right) }+\Delta t\pounds \left( \mathbf{U}%
	^{\left( 2\right) },t^{n}+\frac{1}{2}\Delta t\right) \right) .
\end{eqnarray*}
Here, the operator $\pounds $ represents the DG spatial discretization and $\mathbb{P}_{V^k_h} $ the projection onto the piecewise polynomial space. The stability condition is taken \cite{cockshureview01} as:
\begin{equation*}
	\max\limits_{\forall j}\left( \frac{\max (|\bar{\lambda}^{j}_{1}|,| \bar{\lambda}^j_{2}|)}{%
		\Delta _{j}}\right) \Delta t\leq \frac{1}{2k+1}.
\end{equation*}

\subsection{Slope limiter}  %[{\color{red}this order is better}] }
It is well-known that higher-order space discretizations are prone to generating spurious oscillations \cite{Harten83,bookLev}. Therefore, a
limiting technique is often  needed in a high-order scheme to make the solution stable and oscillation free.
In this study, an usual  slope limiter is applied to the approximate solution after every stage of the Runge-Kutta scheme as done in the RKDG methods \cite{CS98,cockshureview01}.
It is remarked that the initial approximation $\mathbf{U}_{h}(x,0)$ obtained by the projection $\mathbb{P}_{V^k_h} \textbf{U}_0$  is also be limited to suppress possible spurious oscillations.

The  minmod ($mm$) limiter function is defined as
$$mm\left( a_{1},a_{2},a_{3}\right) =\left\{
\begin{array}{ll}
	\sigma \cdot  \min\limits_{1\leq i\leq 3}\left\vert a_{i}\right\vert & \text{if~~} \text{sign}\left( a_{1}\right) =\text{sign}\left( a_{2}\right) =\text{sign}\left(
	a_{3}\right) =\sigma, \\
	0 &  \text{otherwise}. %
\end{array}%
\right. $$
For the higher-order DG solution $\textbf{U}_h$, we first define its linear part $\mathbf{U}_{h}^{1}$ as
\[
\mathbf{U}_{h}^1(x,t)=\mathbf{U}_{j}^{0}   +\mathbf{U}_{j}^{1}\phi_{j}^{1} (x),
\]%
and then define the generalized slope limiter $\Lambda \Pi^k$ acting on the linear part as
\begin{equation*}
	\Lambda \Pi^1  (\mathbf{U}_{h}^{1})=\mathbf{U}_{j}^{0}+mm\left( \textbf{U}^1_j,\gamma (\mathbf{U}_{j}^{0}-\mathbf{U}%
	_{j-1}^{0}),\gamma (\mathbf{U}_{j+1}^{0}-\mathbf{U}_{j}^{0})\right) \frac{%
		(x-x_{j})}{\Delta _{j}/2}.  \label{glimiter}
\end{equation*}
Further, we define two limited interface values as
\begin{equation} \label{half1}
\begin{split}
	\left(\widetilde{\mathbf{U}}_{h}^{1} \right) _{j+1/2}^{-} &= \mathbf{U}
	_{j}^{0}+mm\left( {\mathbf{U}_{j+1/2}^{1-}}   -\mathbf{U}_{j}^{0},~\gamma (
	\mathbf{U}_{j}^{0}-\mathbf{U}_{j-1}^{0}),~\gamma (\mathbf{U}_{j+1}^{0}-
	\mathbf{U}_{j}^{0})\right) ,  \\
	\left(\widetilde{\mathbf{U}}_{h}^{1} \right) _{j+1/2}^{+} &= \mathbf{U}
	_{j}^{0}-mm\left( \mathbf{U}_{j}^{0}-{\mathbf{U}_{j-1/2}^{1+}},~\gamma (
	\mathbf{U}_{j}^{0}-\mathbf{U}_{j-1}^{0}),~\gamma (\mathbf{U}_{j+1}^{0}-
	\mathbf{U}_{j}^{0})\right).
\end{split}
\end{equation}
The limiting procedure of the higher-order DG solution is as follows:

\noindent 1)   By using Eq.~\eqref{half1}, compute $ \left(\widetilde{\mathbf{U}}_{h}^{1}\right) _{j+1/2}^{-}$ and  $ \left (
\widetilde{\mathbf{U}}_{h}^{1} \right) _{j+1/2}^{+}$  respectively.

\noindent 2)   If  $ {\mathbf{U}_{j+1/2}^{1-} }=\left(\widetilde{\mathbf{U}}_{h}^{1}\right) _{j+1/2}^{-}$ and
                ${\mathbf{U}_{j-1/2}^{1+} }=\left(\widetilde{\mathbf{U}}_{h}^{1}\right) _{j+1/2}^{+}$,
		then we do not limit the solution,  i.e.,  $\Lambda\Pi^k \left(\mathbf{U}_{h}|_{I_{j}}\right)= \mathbf{U}_{h}|_{I_{j}}$. 
		
\noindent 3)  If the  condition in  2) is not satisfied,    then $\Lambda \Pi^k \left(\mathbf{U}_{h}|_{I_{j}}\right)=\Lambda \Pi^1 (\mathbf{U}_{h}^{1})$, and high-order degrees of
		freedom are set to zero.

In Eq.~\eqref{half1}, usually $\frac{1}{2}\leq \gamma \leq 1$.  A  fixed value $~\gamma =0.5 $ is used in this study, corresponding to a strict limiter.

%\textbf{However, acceleration and friction terms require an appropriate numerical treatment of the source terms. }

%\textbf{Stability:}

%As we can see that too large of the time step leads to issue. why?

%what happen if the parcel of quantity moves fast enough to skip elements between trimesters?

%we need to ensure
%\[
%CFL=\frac{dt}{dx}\leq C_{\max }=1
%\]%
%for the explicit methods. Implicit methods relax this restriction to permit higher values. In our case $C_{\max }=\max \left( u\pm \sqrt{\beta h}\right)$.

\subsection{Numerical treatments of semi-wet areas and reposing states}\label{sec:3.4}
In many applications, the region covered by the granular materials has a finite extension and is limited by a free boundary which moves with the flow velocity. The region where the avalanche mass has not reached is called a dry area, and the region where the flow depth is significant (we require $h>h_{\rm semi}=10^{-6}$) is called a wet area. In the initial conditions, the flow depth and velocity are zero in the dry area.  However, as the solution advances in time, a numerically  semi-wet area between the dry and wet areas will occur, in which the flow depth is nonzero but less than $h_{\rm semi}$. The semi-wet area often leads to a falsely fast prorogation of the avalanche front. To resolve the issue, a numerical treatment technique is adopted here:  when the computed average depth $\bar{h}_j^{n+1}$ is less than a very small threshold (we take as $h_\epsilon=10^{-10}$), the current cell is regarded as dry area for which the avalanche depth is as computed but the velocity is forced to be zero. When  $h_\epsilon<\bar{h}^{n+1}_j < h_{\rm semi}$, the current cell is regarded as semi-wet area for which the depth is as computed but the velocity is extrapolated in  the following way \cite{liyuan15}
%one of the following three ways;
%when the updated depth is greater than the upper threshold $h_{semi}$, the velocity is computed by  $\displaystyle u=q/h$. The three ways to treat semi-wet cells are:
%1) When the current cell's depth meets $0<h_{i}<h_{semi}=10^{-6}$, the cell is marked as semi-dry cell, the velocity of this cell is obtained as follows \cite{liyuan15}:
\begin{equation*}
\bar{u}_{j} =\left\{
\begin{array}{ll}
0 & \quad \text{if}\quad \bar{h}_{j-1}<h_{\rm semi}\text{ \ and \ }\bar{h}_{j+1}<h_{\rm semi}
 \\
\bar{u}_{j-1}    & \quad \text{else if}\quad \bar{h}_{j-1}\geq \bar{h}_{j+1}, \\
\bar{u}_{j+1}   & \quad \text{else.}%
\end{array}%
\right.
\end{equation*}
The cell-averaged momentum, $q^j_0(t)$, is modified accordingly, and the higher modes of the momentum, $q^j_l(t), l=1,2$, are set to zero.

If the slope of the avalanche surface is less than the friction angle, the materials will stop. The reposing state is implemented only based on the cell-averaged flow variables (the 0th-mode DOFs of the DG solution). Thus, the reposing state of a cell is achieved if the admissible basal friction force is less than the Coulomb friction threshold. The detail can be found in \cite{liyuan15}. For a reposing  cell $M^{\rm stop}_j=0$, otherwise $M^{\rm stop}_j=1$.

\section{Numerical results and discussions}
In this section, numerical results for finite granular mass sliding down an inclined plane and merging continuously into horizontal plane are presented \cite{wang2004}. The computational domain is $\left[0,30\right]$ with cell number $N=256$, time step $\Delta t=0.001$, $\epsilon =\frac{1.85}{30}$ , $\kappa=\frac{\zeta_{0}}{4}$ and degree of the polynomials is $k=2$. The inclined region lies  $x \in [0~~ 17.5)$, and the horizontal run-out region is $x > 21.5$ connected by a smooth  transition zone $x \in [17.5~~ 21.5]$. The granular mass is suddenly released at $t = 0$ from the semi-circular shell with an initial radius of $r_{0}=1.85$ with center located at $x_{0}=4$ in dimensionless length units as shown in  figure \ref{top}.

\subsection{Effects of phenomenological parameters and topography}
\subsubsection{Case I}

In this test case, we set the inclination angle $\zeta_{0}=35^{\circ}$, the internal friction angle $\varphi =30^{\circ}$, and the bed friction angle $\delta =30^{\circ}$. The semi-circular granular mass (avalanche) is released at $t=0$, so that it accelerates in the down-slope direction due to gravity and expands in the inclined region for $t=12$. It can be seen from figure \ref{case1} that the granular materials spread over the inclined zone while part of the granular mass reaches the horizontal run-out zone and begins to deposit because of the basal friction. However, the tail is still in the movement at $t = 12$ to $t = 48$. Finally, at $t = 48$, most of the granular mass is accumulated in the transition zone and at the beginning end of the horizontal zone as shown in figure \ref{case1}.
\begin{figure}[!htbp]
	\centering
	\includegraphics[width=0.9 \textwidth ]{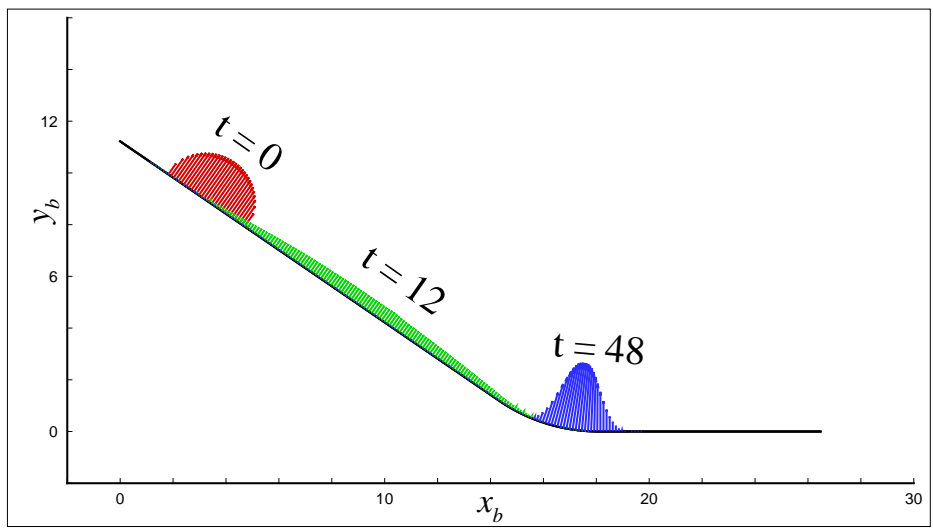}
	\caption{Snapshots of avalanche spatial distribution in case I at different time instants.}\label{case1}
\end{figure}

Figure \ref{fig1}(a) and \ref{fig1}(b) shows  variations of the avalanche depth and velocity at different times  respectively. We can see that the steady state is achieved at $t=48$ marked by $u \doteq 0$ and the
maximum depth  of the avalanche is approximately $2.5$.

\begin{figure*}[!htpb]
	\centering
	\begin{minipage}{0.49\textwidth}
		\includegraphics[width=1.00 \textwidth]{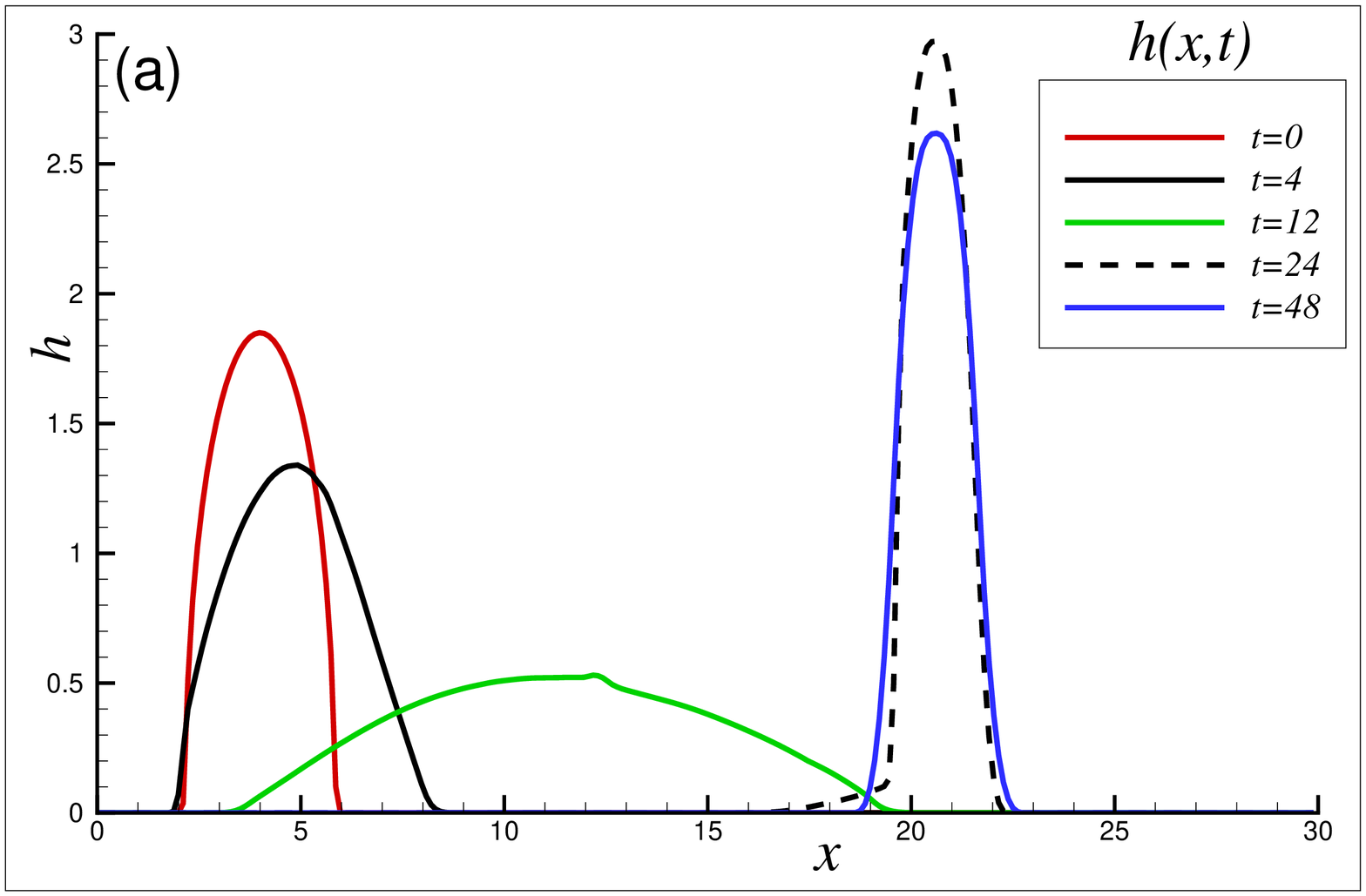}
	\end{minipage}
	\makebox[ 0.00 cm] {}
	\begin{minipage}{0.49\textwidth}
		\includegraphics[width=1.00 \textwidth] {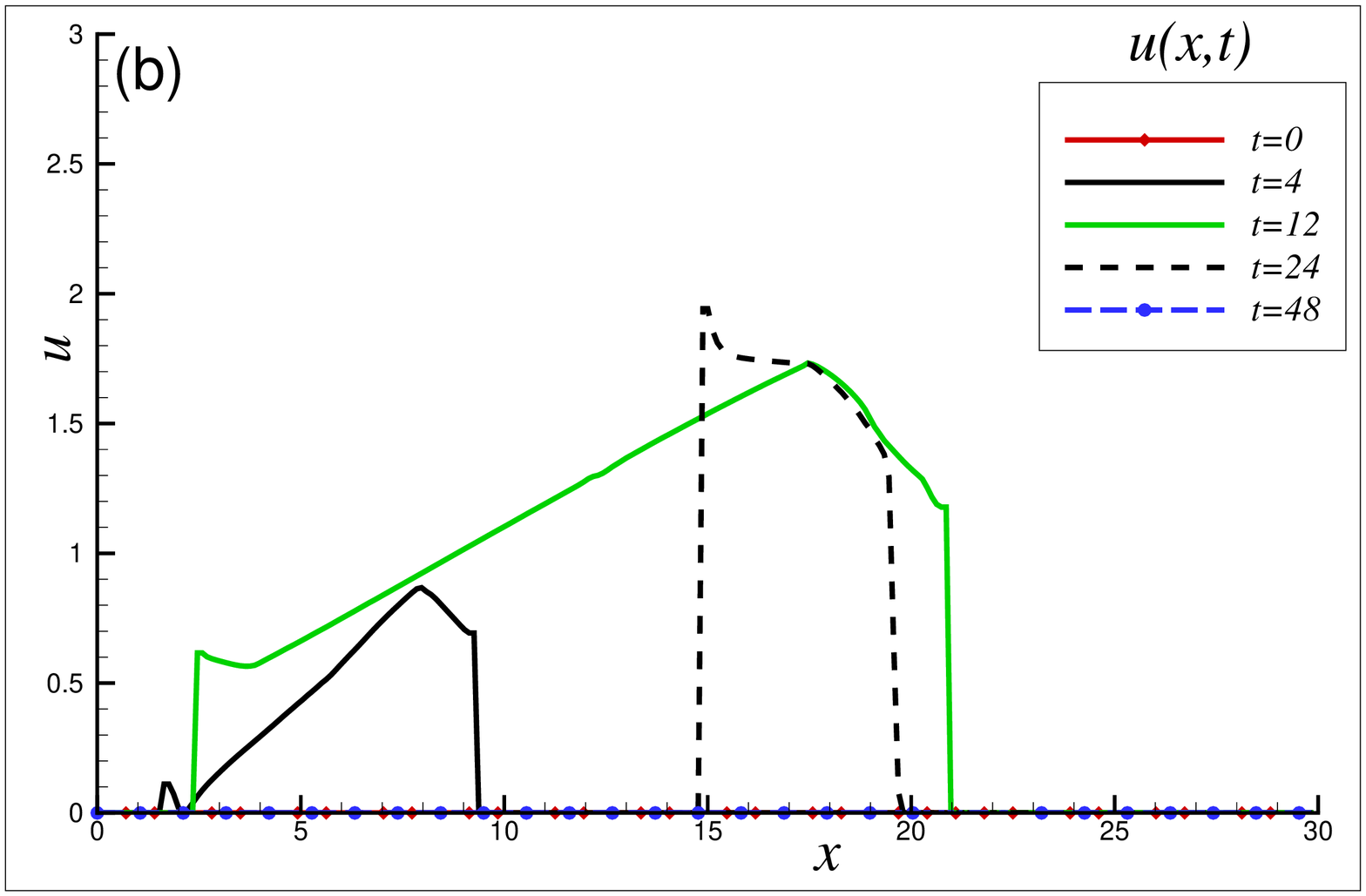}
	\end{minipage}
	\caption{Variations in (a) avalanche depth and (b) velocity of case I at different time instants.}\label{fig1}
\end{figure*}

\subsubsection{Case II}
In this test case, the sensitivity of the granular flow to the bed friction angle  is examined by depicting the evolution of the avalanche with time. We set $\zeta_{0}=35^{\circ}$, $\varphi =30^{\circ}$, but $\delta =23^{\circ}$. From figure \ref{case2} it is observed that with the decrease of the bed friction angle, the granular body becomes more fluidized. Compared with figure \ref{case1},  the run-out zone is larger, and the deposit becomes shallower. This shows that the avalanche is sensitive to the variation of the bed fraction angle $\delta$.

\begin{figure}[!htbp]
	\centering
	\includegraphics[width=0.8 \textwidth ]{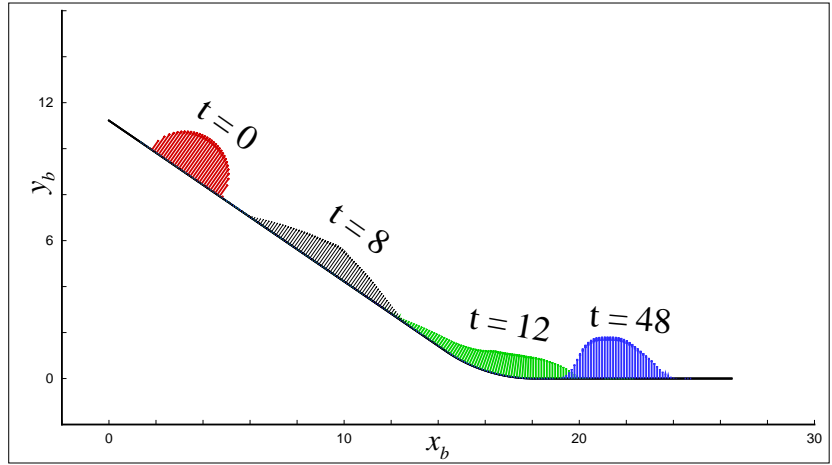}
	\caption{Snapshots of avalanche spatial position for case II at different time instants.}\label{case2}
\end{figure}

Figure~\ref{fig2}(a) illustrates the evolution of the flow depth with time.  The final maximum depth at $t=48$ is $h \simeq  1.8$, smaller than $h\simeq 2.5$ in figure \ref{fig1} of test case I. Further, figure \ref{fig2}(b) shows that the granular materials attain the steady state at time $t=24$ earlier than test case I.

\begin{figure*}[!htb]
	\centering
	\begin{minipage}{0.49 \textwidth}
		\includegraphics[width=1.00 \textwidth]{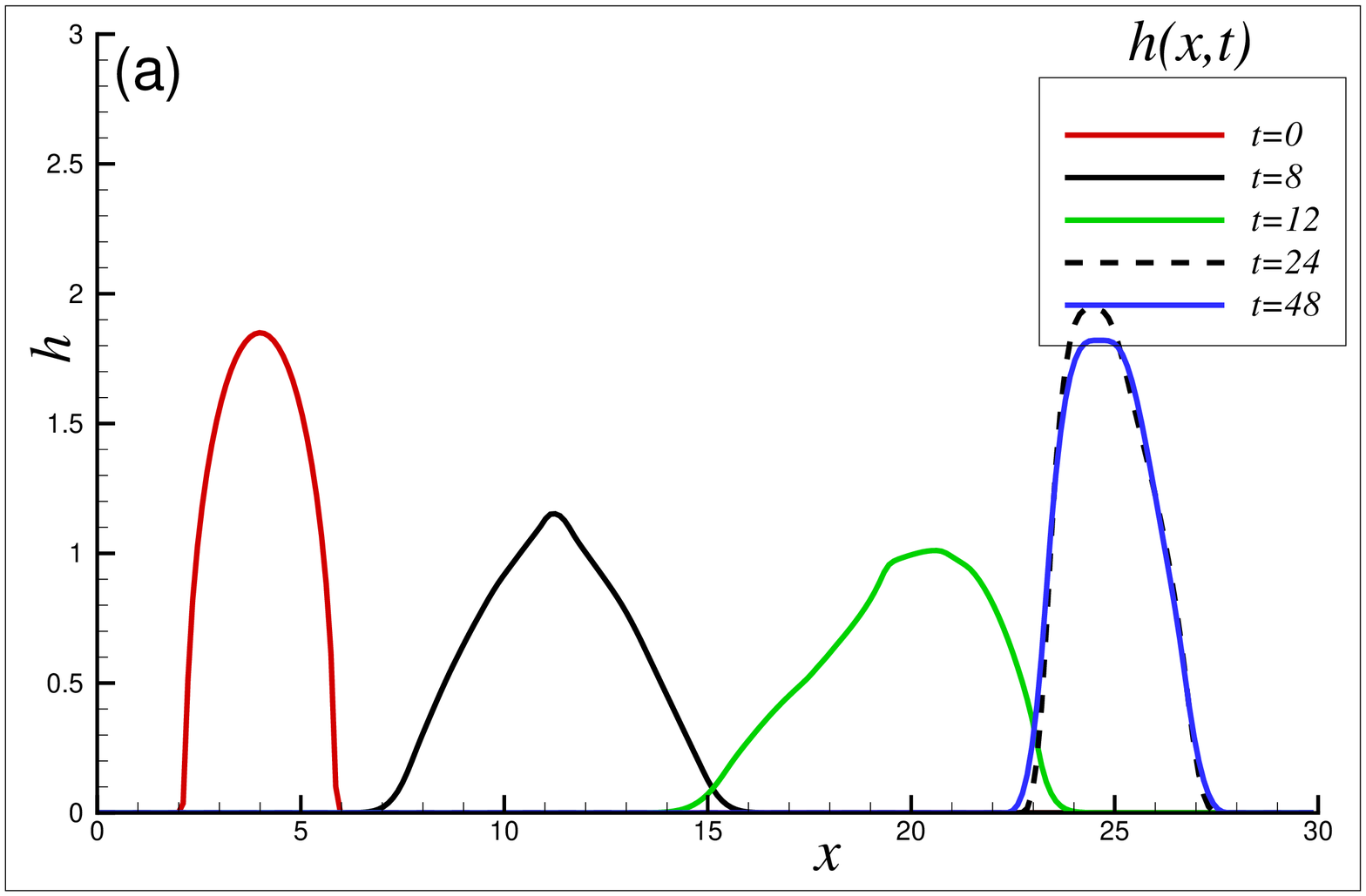}
	\end{minipage}
	\makebox[0.00 cm] {}
	\begin{minipage}{0.49 \textwidth}
		\includegraphics[width=1.00 \textwidth]{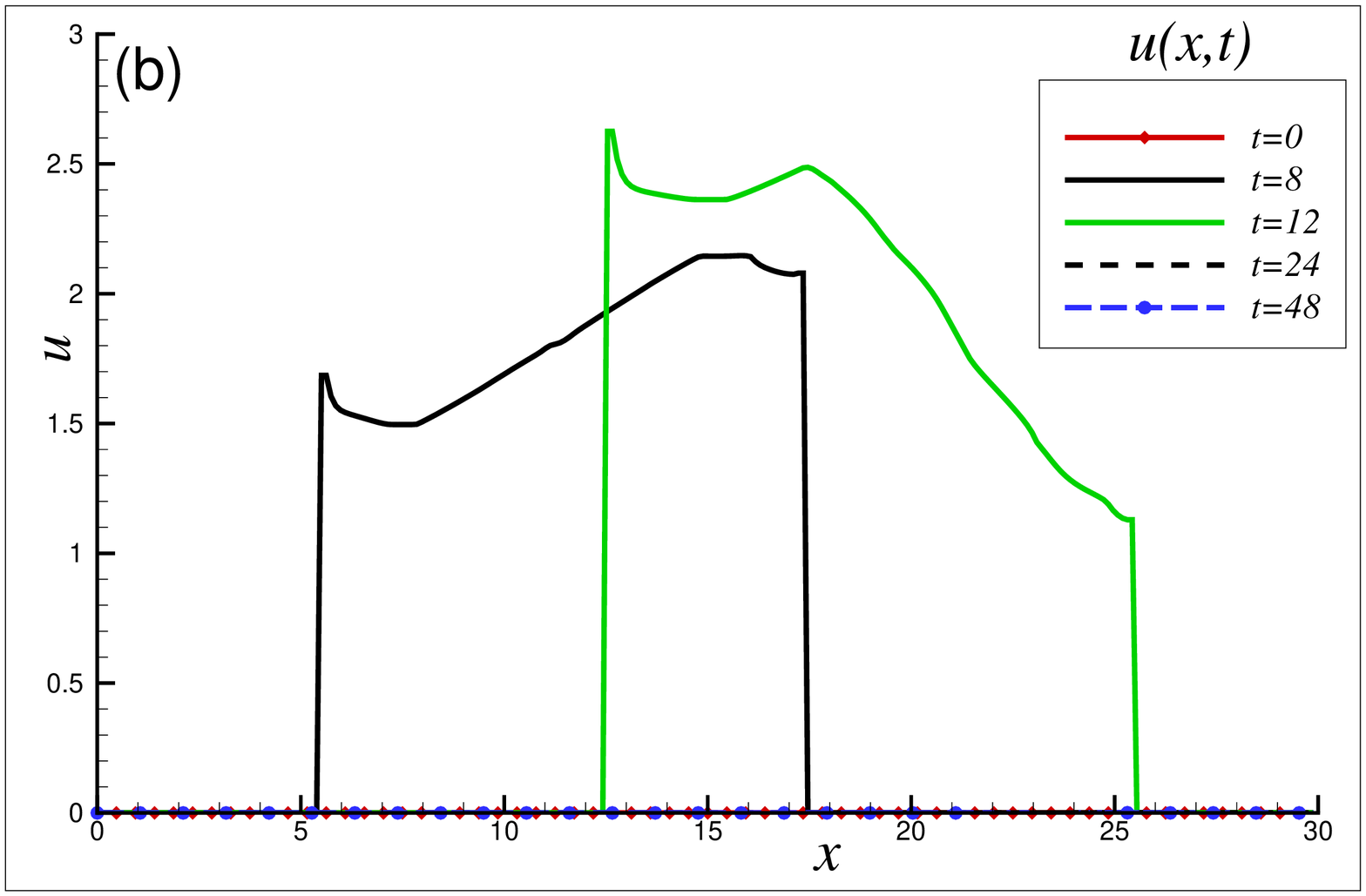}
	\end{minipage}
	\caption{Variations in (a) avalanche height and (b) velocity for case II  at different time instants. }\label{fig2}
\end{figure*}
\subsubsection{Case III}
In this test case, the sensitivity of the granular flow to the internal friction angle $\varphi$  is examined by depicting the evolution of the avalanche with time. We set $\zeta_{0}=35^{\circ}$, $\varphi =37^{\circ}$, and $\delta =30^{\circ}$. In figure \ref{case3}, it is shown that the avalanche is less sensitive to the variation of the internal friction angle $\varphi$ compared with figure \ref{case1}.

%At this stage, a surge wave is created at $t = 12$, which moves a short distance upward as can clearly be seen by comparing the humps for $t = 12$ to $t = 48$.

\begin{figure}[htb!]
	\centering
	\includegraphics[width=0.8 \textwidth  ]{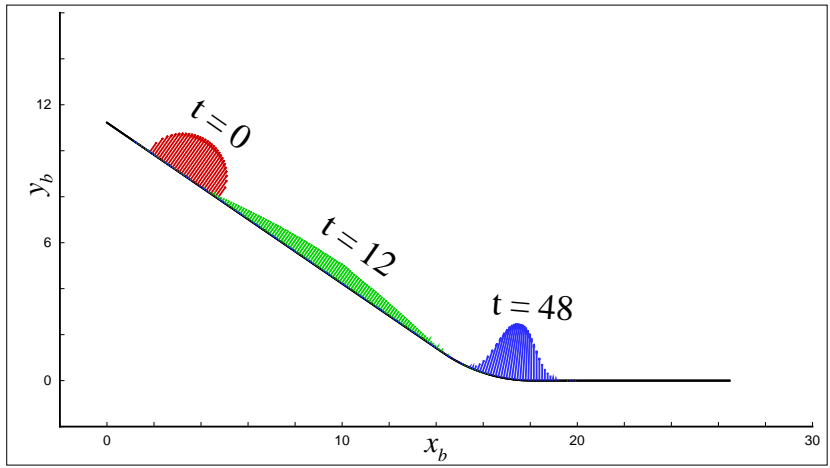}
	\caption{Snapshots of avalanche spatial position for case III at different time instants. }\label{case3}
\end{figure}

 With the increase of the internal friction angle, no much variation in the run-out zone or the avalanche depth  is observed as shown in figure \ref{fig3}(a) compared with a  smaller internal friction angle given in figure \ref{fig1}(a) of test case I.  The only difference between  test cases I and  III is the  avalanche depth becomes a little smaller  with the increasing of  internal friction angle $\varphi$.

The variation of avalanche velocity is illustrated in figure~\ref{fig3}(b). We can see that the granular material attains the steady state a bit earlier, for example, compare the velocity at  $t=24$ with that in figure \ref{fig1}(b) for $\varphi=30^{\circ}$.

\begin{figure*}[htbp!]
	\centering
	\begin{minipage}{0.49 \textwidth}
		\includegraphics[width=1.00 \textwidth]{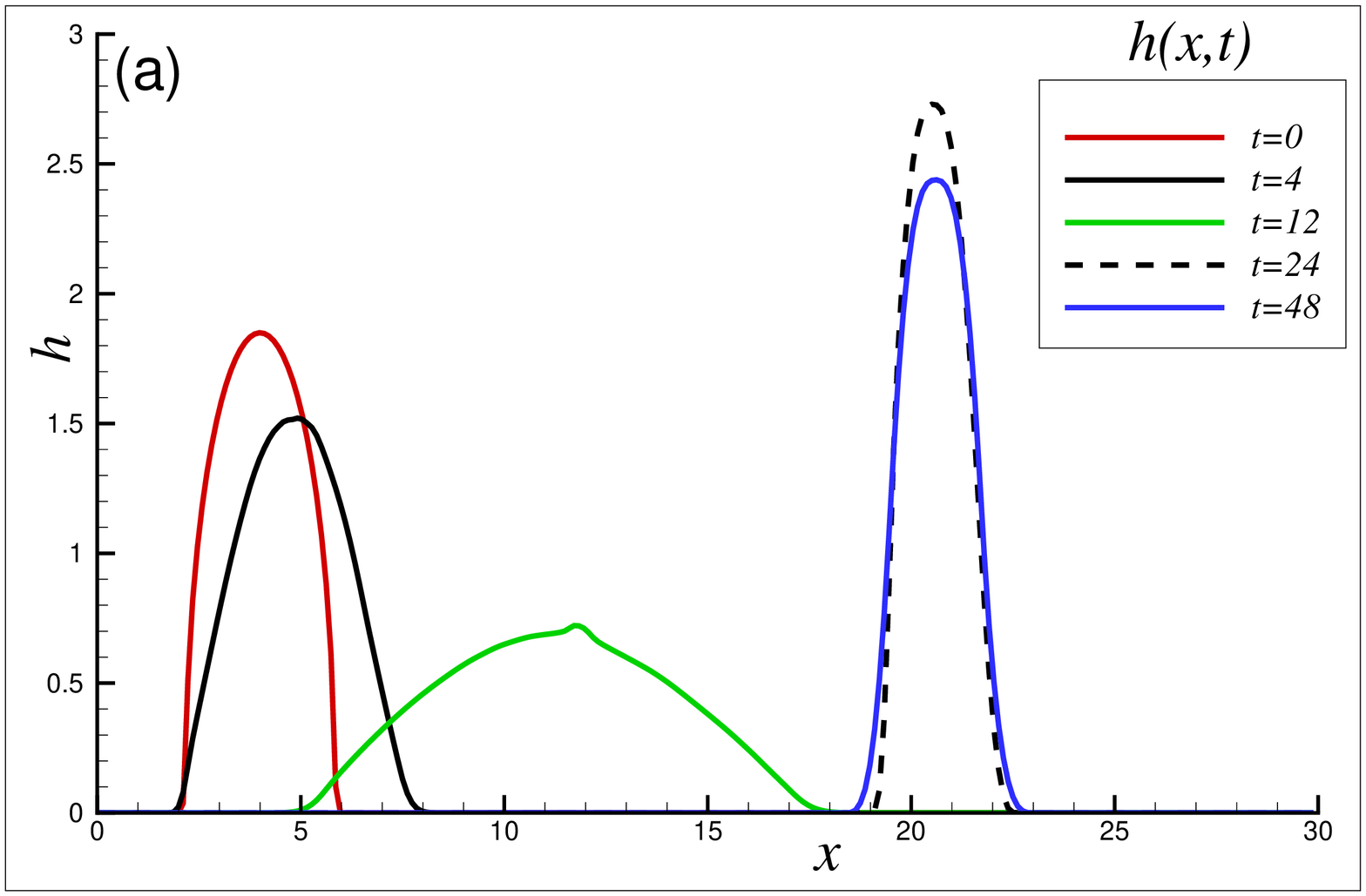}
	\end{minipage}
	\makebox[0.00 cm] {}
	\begin{minipage}{0.49 \textwidth}
		\includegraphics[width=1.00 \textwidth]{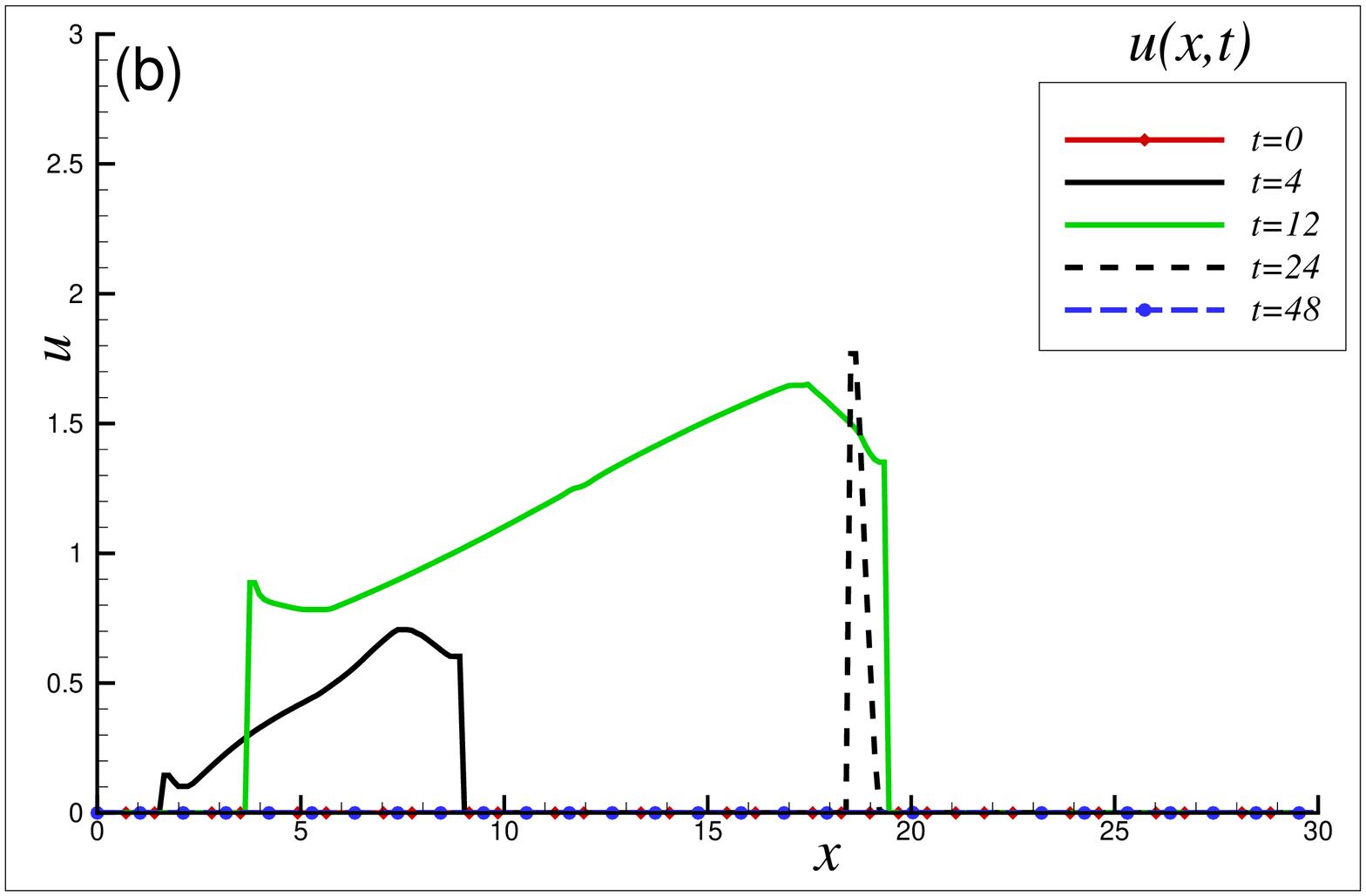}
	\end{minipage}
	\caption{Variation in (a) avalanche height and (b) velocity for case III at different time instants.}  \label{fig3}
\end{figure*}

\subsection{Case IV}

In this test case, we examine the effect of the inclination angle and set $\zeta_{0}=40^{\circ}$, $\varphi =30^{\circ}$, and $\delta =30^{\circ}$. It is observed that the avalanche spreads at a faster speed and goes  a longer distance in the run-out zone as shown in figure \ref{case4}.

\begin{figure}[htb!]
	\centering
	\includegraphics[width=0.8\textwidth]{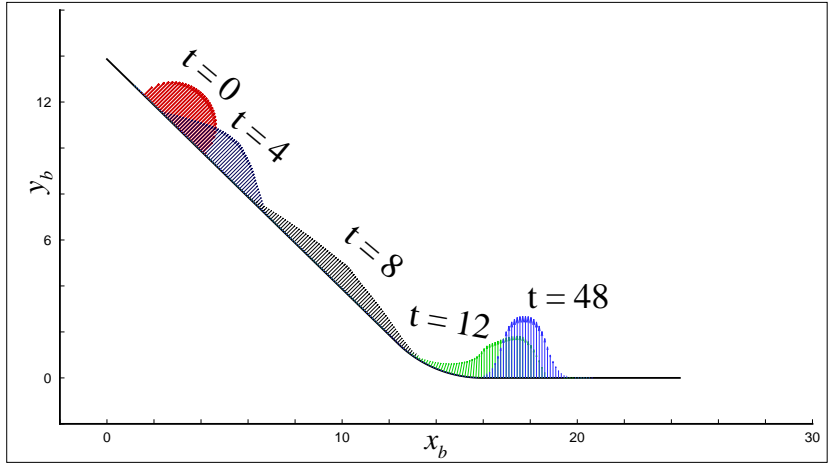}
	\caption{Snapshots of avalanche spatial position  for case IV at different time instants.}\label{case4}
\end{figure}

The variations of avalanche height and velocity are illustrated in figure \ref{fig4}(a) and figure \ref{fig4}(b) respectively. It is noted that the avalanche has a larger height in compassion with test case I, and it reaches earlier to the steady state at $t=24$ as shown in figure \ref{fig4}(b).

\begin{figure*}[htb!]
	\centering
	\begin{minipage}{0.49 \textwidth}
	\includegraphics[width=1.0 \textwidth]{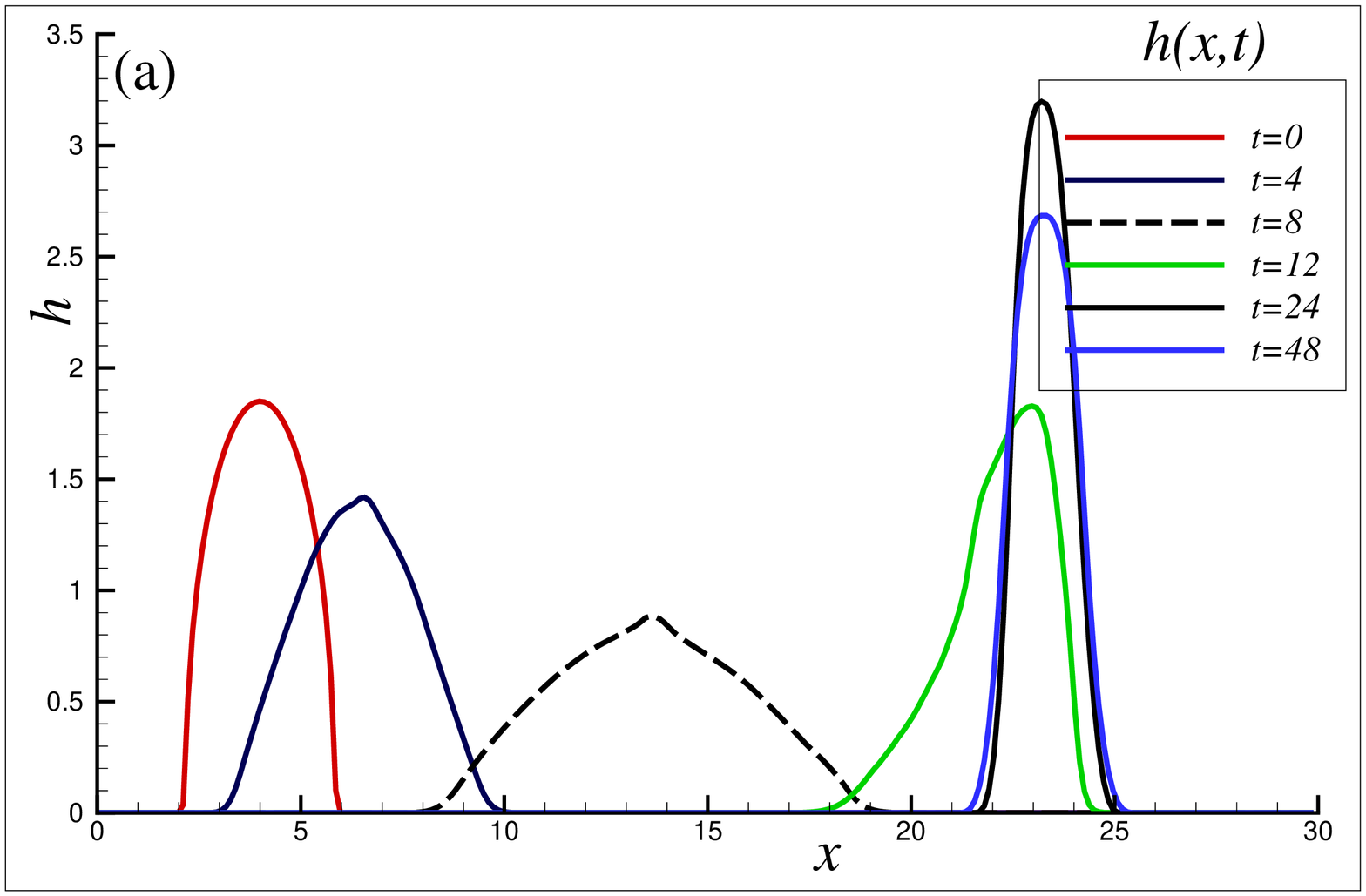}
	\end{minipage}
	\makebox[0.00 cm] {}
	\begin{minipage}{0.49 \textwidth}
		\includegraphics[width=1.0 \textwidth] {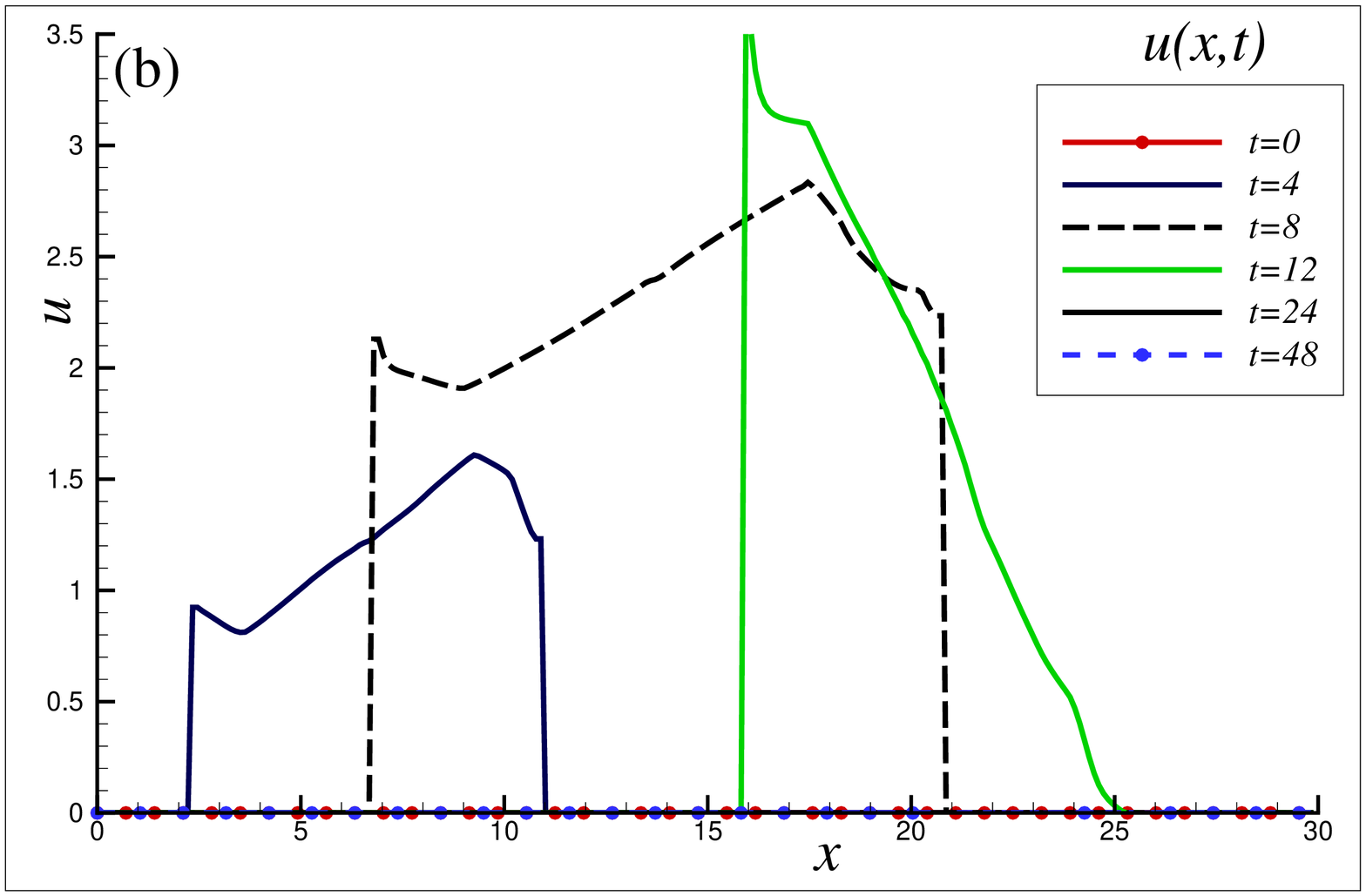}
	\end{minipage}
	\caption{Variations in (a) avalanche height and (b) velocity for case IV at different time instants.}\label{fig4}
\end{figure*}
\section{Conclusion}
This work is an attempt to apply a Runge-Kutta discontinuous Galerkin (RKDG)  method to numerically solve the Savage-Hutter equations that describes the motion of granular avalanche on a curved slope. A third-order space and time accuracy RKDG method is used to solve the one-dimensional SH equations  formulated in a reference surface. Details of the DG discretizations, limiter, and numerical treatments of dry-wet fronts and
reposing states of granular materials are given.

The parameter investigations show that the avalanche flow is more sensitive to variations of the bed friction and inclination angles. It is observed that by decreasing the bed friction angle, the avalanche accelerates more rapidly and attains higher speeds on the inclined bed. It also becomes more shallows, so maximum run-out distance is therefore increased and the avalanche height decreased. In the case of an increasing inclination angle, the avalanche also moves faster on the inclined bed and comes to rest quickly but the height of the avalanche increases. The internal friction angle has a small effect on the movement of the avalanche. By increasing the internal friction angle, the deposit depth decreases slightly.

Extension of the present method to two-dimensional problems having variation in topography in the cross-slope direction is our ongoing work.

\section{Acknowledgments} This work was carried out while A. Shah was visiting ICMSEC(AMSS), Chinese Academy of Sciences, Beijing China as a PIFI visiting scientist. L.~Yuan was supported by the state key program for developing basic sciences ($2010CB731505$)  and the Natural Science Foundation of China ($11261160486$, $91641107$, $91852116$). The work of A. Shah was supported by HEC under NRPU No. $7781$.

\bibliographystyle{plain}
\bibliography{mybibfile}

%\begin{thebibliography}{00}

%
%\bibitem{Reed1973} W. H. Reed and T. R. Hill, Triangular mesh methods for the neutron transport equation, Los Alamos Scientific Lab, N. Mex.(USA), (1973) LA-UR 73-479.%
%

%\end{thebibliography}

\section{Appendix}
The following algebraic relations are used for coordinate transformation
from the computational domain $(x,h)$ to the {horizontal-vertical} domain
with $x=(x_{b},y_{b})$, and $h=(h_{x_{b}},h_{y_{b}})$. Let $r=\frac{4}{\zeta
	_{0}}$ and define 
\begin{eqnarray*}
	&&\left. x_{1}=17.5\cos \left( \zeta _{0}\right) ;~~~y_{1}=r\left( 1-\cos
	\left( \zeta _{0}\right) \right) ,\right.  \\
	&&\left. x_{2}=x_{1}+r\sin \left( \zeta _{0}\right) ;~~~y_{2}=0.\right. 
\end{eqnarray*}%
{We assume that the projected }$x_{b}{-}${axis is aligned with the
	horizontal part of the reference surface, then }the transformation is
follows:\\
If
\begin{equation*}
\text{ }x\geq 21.5\text{ then }x_{b}=x_{2}+(x-21.5),\text{ \ }y_{b}=0,%
\text{  }h_{x_{b}}=0,\text{~and }h_{y_{b}}=h,
\end{equation*}%
else if
\begin{eqnarray*}
	&&\left. \text{ }x\geq 17.5\text{ and }x<21.5\text{ then }\zeta =\left(
	21.5-x\right) /r,\text{ \ }x_{b}=x_{2}-r\sin (\zeta ),\right.  \\
	&&\left. y_{b}=r(1-\cos (\zeta )),\text{ }h_{x_{b}}=h\sin (\zeta ),\text{
		and }h_{y_{b}}=h\cos (\zeta ),\right. 
\end{eqnarray*}%
else if
\begin{eqnarray*}
	&&\left. \text{ }x<17.5\text{ then }x_{b}=x\cos \left( \zeta _{0}\right) ,%
	\text{ \ }y_{b}=y_{1}+\left( 17.5-x\right) \sin \zeta _{0},\right.  \\
	&&\left. h_{x_{b}}=h\sin \left( \zeta _{0}\right) ,\text{ and }%
	h_{y_{b}}=h\cos \left( \zeta _{0}\right) ,\right. 
\end{eqnarray*}
\text{end if.}
\end{document}